\documentclass[reqno,12pt, a4paper]{amsart}
\usepackage[foot]{amsaddr}

\usepackage{mathrsfs}
\usepackage{amssymb}
\usepackage{amsfonts}
\usepackage{latexsym}
\usepackage{amsthm}

\usepackage{graphicx}
\def\lb{\label}

\newcommand{\er}[1]{\textrm{(\ref{#1})}}

\newcommand{\HOX}[1]{\marginpar{\footnotesize #1}}

\begin{document}


\renewcommand{\theequation}{\arabic{section}.\arabic{equation}}
\theoremstyle{plain}
\newtheorem{theorem}{\bf Theorem}[section]
\newtheorem{lemma}[theorem]{\bf Lemma}
\newtheorem{corollary}[theorem]{\bf Corollary}
\newtheorem{proposition}[theorem]{\bf Proposition}
\newtheorem{definition}[theorem]{\bf Definition}
\newtheorem{remark}[theorem]{\it Remark}

\def\a{\alpha}  \def\cA{{\mathcal A}}     \def\bA{{\bf A}}  \def\mA{{\mathscr A}}
\def\b{\beta}   \def\cB{{\mathcal B}}     \def\bB{{\bf B}}  \def\mB{{\mathscr B}}
\def\g{\gamma}  \def\cC{{\mathcal C}}     \def\bC{{\bf C}}  \def\mC{{\mathscr C}}
\def\G{\Gamma}  \def\cD{{\mathcal D}}     \def\bD{{\bf D}}  \def\mD{{\mathscr D}}
\def\d{\delta}  \def\cE{{\mathcal E}}     \def\bE{{\bf E}}  \def\mE{{\mathscr E}}
\def\D{\Delta}  \def\cF{{\mathcal F}}     \def\bF{{\bf F}}  \def\mF{{\mathscr F}}
\def\c{\chi}    \def\cG{{\mathcal G}}     \def\bG{{\bf G}}  \def\mG{{\mathscr G}}
\def\z{\zeta}   \def\cH{{\mathcal H}}     \def\bH{{\bf H}}  \def\mH{{\mathscr H}}
\def\e{\eta}    \def\cI{{\mathcal I}}     \def\bI{{\bf I}}  \def\mI{{\mathscr I}}
\def\p{\psi}    \def\cJ{{\mathcal J}}     \def\bJ{{\bf J}}  \def\mJ{{\mathscr J}}
\def\vT{\Theta} \def\cK{{\mathcal K}}     \def\bK{{\bf K}}  \def\mK{{\mathscr K}}
\def\k{\kappa}  \def\cL{{\mathcal L}}     \def\bL{{\bf L}}  \def\mL{{\mathscr L}}
\def\l{\lambda} \def\cM{{\mathcal M}}     \def\bM{{\bf M}}  \def\mM{{\mathscr M}}
\def\L{\Lambda} \def\cN{{\mathcal N}}     \def\bN{{\bf N}}  \def\mN{{\mathscr N}}
\def\m{\mu}     \def\cO{{\mathcal O}}     \def\bO{{\bf O}}  \def\mO{{\mathscr O}}
\def\n{\nu}     \def\cP{{\mathcal P}}     \def\bP{{\bf P}}  \def\mP{{\mathscr P}}
\def\r{\varrho} \def\cQ{{\mathcal Q}}     \def\bQ{{\bf Q}}  \def\mQ{{\mathscr Q}}
\def\s{\sigma}  \def\cR{{\mathcal R}}     \def\bR{{\bf R}}  \def\mR{{\mathscr R}}
\def\S{\Sigma}  \def\cS{{\mathcal S}}     \def\bS{{\bf S}}  \def\mS{{\mathscr S}}
\def\t{\tau}    \def\cT{{\mathcal T}}     \def\bT{{\bf T}}  \def\mT{{\mathscr T}}
\def\f{\phi}    \def\cU{{\mathcal U}}     \def\bU{{\bf U}}  \def\mU{{\mathscr U}}
\def\F{\Phi}    \def\cV{{\mathcal V}}     \def\bV{{\bf V}}  \def\mV{{\mathscr V}}
\def\P{\Psi}    \def\cW{{\mathcal W}}     \def\bW{{\bf W}}  \def\mW{{\mathscr W}}
\def\o{\omega}  \def\cX{{\mathcal X}}     \def\bX{{\bf X}}  \def\mX{{\mathscr X}}
\def\x{\xi}     \def\cY{{\mathcal Y}}     \def\bY{{\bf Y}}  \def\mY{{\mathscr Y}}
\def\X{\Xi}     \def\cZ{{\mathcal Z}}     \def\bZ{{\bf Z}}  \def\mZ{{\mathscr Z}}
\def\O{\Omega}

\newcommand{\mc}{\mathscr {c}}

\newcommand{\gA}{\mathfrak{A}}          \newcommand{\ga}{\mathfrak{a}}
\newcommand{\gB}{\mathfrak{B}}          \newcommand{\gb}{\mathfrak{b}}
\newcommand{\gC}{\mathfrak{C}}          \newcommand{\gc}{\mathfrak{c}}
\newcommand{\gD}{\mathfrak{D}}          \newcommand{\gd}{\mathfrak{d}}
\newcommand{\gE}{\mathfrak{E}}
\newcommand{\gF}{\mathfrak{F}}           \newcommand{\gf}{\mathfrak{f}}
\newcommand{\gG}{\mathfrak{G}}           
\newcommand{\gH}{\mathfrak{H}}           \newcommand{\gh}{\mathfrak{h}}
\newcommand{\gI}{\mathfrak{I}}           \newcommand{\gi}{\mathfrak{i}}
\newcommand{\gJ}{\mathfrak{J}}           \newcommand{\gj}{\mathfrak{j}}
\newcommand{\gK}{\mathfrak{K}}            \newcommand{\gk}{\mathfrak{k}}
\newcommand{\gL}{\mathfrak{L}}            \newcommand{\gl}{\mathfrak{l}}
\newcommand{\gM}{\mathfrak{M}}            \newcommand{\gm}{\mathfrak{m}}
\newcommand{\gN}{\mathfrak{N}}            \newcommand{\gn}{\mathfrak{n}}
\newcommand{\gO}{\mathfrak{O}}
\newcommand{\gP}{\mathfrak{P}}             \newcommand{\gp}{\mathfrak{p}}
\newcommand{\gQ}{\mathfrak{Q}}             \newcommand{\gq}{\mathfrak{q}}
\newcommand{\gR}{\mathfrak{R}}             \newcommand{\gr}{\mathfrak{r}}
\newcommand{\gS}{\mathfrak{S}}              \newcommand{\gs}{\mathfrak{s}}
\newcommand{\gT}{\mathfrak{T}}             \newcommand{\gt}{\mathfrak{t}}
\newcommand{\gU}{\mathfrak{U}}             \newcommand{\gu}{\mathfrak{u}}
\newcommand{\gV}{\mathfrak{V}}             \newcommand{\gv}{\mathfrak{v}}
\newcommand{\gW}{\mathfrak{W}}             \newcommand{\gw}{\mathfrak{w}}
\newcommand{\gX}{\mathfrak{X}}               \newcommand{\gx}{\mathfrak{x}}
\newcommand{\gY}{\mathfrak{Y}}              \newcommand{\gy}{\mathfrak{y}}
\newcommand{\gZ}{\mathfrak{Z}}             \newcommand{\gz}{\mathfrak{z}}

\def\ve{\varepsilon}   \def\vt{\vartheta}    \def\vp{\varphi} \def\vk{\varkappa}

\def\A{{\mathbb A}} \def\B{{\mathbb B}} \def\C{{\mathbb C}}\def\dD{{\mathbb D}}
\def\E{{\mathbb E}} \def\dF{{\mathbb F}} \def\dG{{\mathbb G}} \def\H{{\mathbb H}}
\def\I{{\mathbb I}} \def\J{{\mathbb J}} \def\K{{\mathbb K}} \def\dL{{\mathbb L}}
\def\M{{\mathbb M}} \def\N{{\mathbb N}} \def\O{{\mathbb O}} \def\dP{{\mathbb P}}
 \def\R{{\mathbb R}} \def\dQ{{\mathbb Q}}
\def\dS{{\mathbb S}} \def\T{{\mathbb T}} \def\U{{\mathbb U}} \def\V{{\mathbb V}}
\def\W{{\mathbb W}} \def\X{{\mathbb X}} \def\Y{{\mathbb Y}} \def\Z{{\mathbb Z}}

\newcommand{\1}{\mathbbm 1}
\newcommand{\dd}    {\, \mathrm d}



\def\la{\leftarrow}              \def\ra{\rightarrow} \def\Ra{\Rightarrow}
\def\ua{\uparrow}                \def\da{\downarrow}
\def\lra{\leftrightarrow}        \def\Lra{\Leftrightarrow}


\def\lt{\biggl}                  \def\rt{\biggr}
\def\ol{\overline}               \def\wt{\widetilde}
\def\no{\noindent}


\let\ge\geqslant                 \let\le\leqslant
\def\lan{\langle}                \def\ran{\rangle}
\def\/{\over}                    \def\iy{\infty}
\def\sm{\setminus}               \def\es{\emptyset}
\def\ss{\subset}                 \def\ts{\times}
\def\pa{\partial}                \def\os{\oplus}
\def\om{\ominus}                 \def\ev{\equiv}
\def\iint{\int\!\!\!\int}        \def\iintt{\mathop{\int\!\!\int\!\!\dots\!\!\int}\limits}
\def\el2{\ell^{\,2}}             \def\1{1\!\!1}
\def\wh{\widehat}

\def\sh{\mathop{\mathrm{sh}}\nolimits}
\def\ch{\mathop{\mathrm{ch}}\nolimits}
\def\all{\mathop{\mathrm{all}}\nolimits}
\def\where{\mathop{\mathrm{where}}\nolimits}
\def\as{\mathop{\mathrm{as}}\nolimits}
\def\Area{\mathop{\mathrm{Area}}\nolimits}
\def\arg{\mathop{\mathrm{arg}}\nolimits}
\def\const{\mathop{\mathrm{const}}\nolimits}
\def\det{\mathop{\mathrm{det}}\nolimits}
\def\diag{\mathop{\mathrm{diag}}\nolimits}
\def\diam{\mathop{\mathrm{diam}}\nolimits}
\def\dim{\mathop{\mathrm{dim}}\nolimits}
\def\dist{\mathop{\mathrm{dist}}\nolimits}
\def\Im{\mathop{\mathrm{Im}}\nolimits}
\def\Iso{\mathop{\mathrm{Iso}}\nolimits}
\def\Ker{\mathop{\mathrm{Ker}}\nolimits}
\def\Lip{\mathop{\mathrm{Lip}}\nolimits}
\def\rank{\mathop{\mathrm{rank}}\limits}
\def\Ran{\mathop{\mathrm{Ran}}\nolimits}
\def\Re{\mathop{\mathrm{Re}}\nolimits}
\def\Res{\mathop{\mathrm{Res}}\nolimits}
\def\res{\mathop{\mathrm{res}}\limits}
\def\sign{\mathop{\mathrm{sign}}\nolimits}
\def\span{\mathop{\mathrm{span}}\nolimits}
\def\supp{\mathop{\mathrm{supp}}\nolimits}
\def\Tr{\mathop{\mathrm{Tr}}\nolimits}
\def\BBox{\hspace{1mm}\vrule height6pt width5.5pt depth0pt \hspace{6pt}}


\newcommand\nh[2]{\widehat{#1}\vphantom{#1}^{(#2)}}
\def\dia{\diamond}

\def\Oplus{\bigoplus\nolimits}




\def\qqq{\qquad}
\def\qq{\quad}
\let\ge\geqslant
\let\le\leqslant
\let\geq\geqslant
\let\leq\leqslant
\newcommand{\ca}{\begin{cases}}
\newcommand{\ac}{\end{cases}}
\newcommand{\ma}{\begin{pmatrix}}
\newcommand{\am}{\end{pmatrix}}
\renewcommand{\[}{\begin{equation}}
\renewcommand{\]}{\end{equation}}
\def\bu{\bullet}

\baselineskip 14pt

\title[{Inverse spectral theory for perturbed torus }]
{Inverse spectral theory for perturbed torus}

\date{\today}

\author[Hiroshi Isozaki]{Hiroshi Isozaki}
\address{Professor Emeritus,
University of Tsukuba, Tsukuba, 305-8571, Japan, isozakih@math.tsukuba.ac.jp}
\author[Evgeny L. Korotyaev]{Evgeny L. Korotyaev}
\address{
Department of Math. Analysis, Saint-Petersburg State University,
Universitetskaya nab. 7/9, St. Petersburg, 199034, Russia, \
korotyaev@gmail.com, \ e.korotyaev@spbu.ru}


\begin{abstract}
\no We consider an inverse problem for Laplacians on rotationally
symmetric  manifolds, which are finite  for the transversal
direction and periodic with respect to the axis of the manifold,
i.e.,  Laplacians on tori.  We construct an infinite dimensional
analytic isomorphism between the space of profiles  (the radius of
the rotation) of the torus and the spectral data as well as the
stability estimates: those for the spectral data in terms of the
profile and conversely, for the profile in term of the spectral
data.

\end{abstract}

\maketitle

\section {Introduction}
\setcounter{equation}{0}

\subsection {Geometry.}
In this paper we study  inverse problems on the manifold $M = \dS^1
\times Y$,
 where $Y$ is a compact  Riemannian manifold of dimension $m$ with
 or without boundary equipped with the metric $g_Y$, which gives the warped product metric
\[
\lb{1} g_M = (dx)^2 + r^2(x)g_Y
\]
on $M$. Here, we identify $\dS^1$ with ${\T}^1 = {\R}^1/{\Z} =
[0,1]$ (see Fig. \ref{Fig1}).  We are mainly interested in the case
in which $Y$ is diffeomorphic to (a part of)  $\dS^1$. Therefore, we
call $M$ a perturbed torus, and Y a transversal manifold. The
Laplacian $\D_M$ on $M$ has the form
$$
- \D_M=-{1\/r(x)^m}\pa_x\Big( r(x)^m \pa_x\Big)-{\D_{Y}\/r(x)^2},
$$
where $\Delta_Y$ is the Laplacian on $Y$, which has a discrete
spectrum $0\le E_0<E_1\le E_2\le  . . .$ and the associated complete orthonormal
 system of eigenfunctions ${\P_\n, \n\ge
0}$, in $L^2(Y )$.
Introduce the Hilbert spaces
$$
\mL_\n^2=\rt\{f(x)\P_\n(y) \, ; \, (x,y)\in \dS^1 \ts Y, \
\int_0^1|f(x)|^2r^m(x)dx < \infty\rt\},\ \n\ge 0.
$$
Then, the Laplacian on $(M, g)$ acting on
$L^2(M)\backsimeq\os_{\n\ge 0} \mL_\n^2$
 is unitarily equivalent to a
direct sum of one-dimensional  operators $\D_\n$, namely,
\[
\lb{3} \D_{M}\backsimeq\os_{\n=0}^{\iy}\ \D_\n.
\nonumber
\]
Here the operator $\D_\n$ acts in the space $L^2([0,1],r^m(x)dx)$
and is given by
\[
\begin{aligned}
\lb{3x}
-\D_\n =-{1\/r(x)^m} \pa_x \Big(r(x)^m \pa_x\Big) +{E_\n\/r(x)^2}.
\end{aligned}
\nonumber
\]


\begin{figure}
\tiny
\unitlength 0.7mm 
\linethickness{0.4pt}
\ifx\plotpoint\undefined\newsavebox{\plotpoint}\fi 
\begin{picture}(188.75,45.5)(0,0)
\qbezier(22.75,27.5)(42.375,15.375)(73.5,26.75)
\qbezier(26.25,25.5)(47.625,32.75)(67.5,25)
\qbezier(3.25,25.25)(3.125,9.25)(41.5,5.25)
\qbezier(3.25,25.25)(3.125,41.25)(41.5,45.25)
\qbezier(41.5,45.25)(89.25,44.5)(91,24.5)
\qbezier(41.5,5.25)(89.25,6)(91,24.5)
\qbezier(113.25,10.75)(102,15)(100.75,25.25)
\qbezier(100.75,25.25)(100.875,41)(139.25,45.25)
\qbezier(139.25,45.25)(187,44.25)(188.75,24.25)
\qbezier(125.75,28.75)(121.625,28.5)(120,29.25)
\qbezier[10](131.25,20.5)(118,24.375)(128.75,26.75)
\qbezier(164.5,24.25)(144.625,31.5)(128.75,26.75)
\qbezier(113.5,10.75)(111,24.375)(120.5,28.0)
\qbezier(131.25,20.25)(130.125,31.5)(120.5,28.0)
\multiput(113.5,11)(.0654545455,.0336363636){275}{\line(1,0){.0654545455}}
\qbezier(164.75,24.5)(164.875,35.75)(176.75,36.25)
\qbezier(176.75,36.25)(188,36.125)(188.25,24.5)
\put(164.75,24.5){\line(1,0){23.75}}
\end{picture}
\caption{\footnotesize The torus} \lb{Fig1}
\end{figure}

\subsection{Example}
Our motivating example is the
 following perturbed torus in ${\R}^3$:
\begin{equation}
x = \left(a + R(\theta)\cos\theta\right)\cos\phi, \quad
y = \left(a + R(\theta)\cos\theta\right)\sin\phi, \quad z
= a + R(\theta)\sin\theta
\nonumber
\end{equation}
where
\[
a > 0, \qq R(\theta) \in C^{2}(\dS^1), \qq 0 < R(\theta) < a, \qqq
\theta, \phi \in [0,2\pi]. \nonumber
\]
The induced metric on this surface is
\begin{equation}
ds^2 =\left(R'(\theta)^2 + R(\theta)^2\right)(d\theta)^2 +
\left(a + R(\theta)\cos\theta\right)^2(d\phi)^2.
\nonumber
\end{equation}
We put
$$
b = \int_0^{2\pi}\sqrt{R'(\theta)^2 + R(\theta)^2}d\theta,
$$
and make the change of variable $t=t(\theta)$ by
$$
\frac{d t}{d\theta} = \sqrt{R'(\theta)^2 + R(\theta)^2}.
$$
Then we have
$$
ds^2 = (d t)^2 + r(t)^2(d\phi)^2, \quad r(t) = a + R(\theta(t))\cos\theta(t),
$$
where $0 \leq \phi\leq 2\pi, 0 \leq t \leq b$. We also have
$$
|r'(t)| \leq 1,
$$
since
$$
r'(t) = \frac{d
\left(R(\theta)\cos\theta\right)}{d\theta}\frac{d\theta}{d t}
 = \frac{R'(\theta)\cos\theta - R(\theta)\sin\theta}{\sqrt{R'(\theta)^2
  + R(\theta)^2}}.
$$
Putting $\t = \frac{t}{b}$, we can rewrite it as
$$
ds^2 =b^2\rt((d\t)^2 + h(\t)^2(d\phi)^2\rt),
 \quad h(\t) = \frac{r(t)}{b},
$$
$$
|h'(\t)| = |r'(t)| \leq 1, \qq \forall \ t\in [0, b].
$$

\subsection{Problem}
This example leads us to the manifold $M= \dS^1\times Y$ equipped
with metric
\begin{equation}
g_M= b^2\left((dx)^2 + r(x)^2g_Y\right),
\end{equation}
where $b > 0$ is a constant.
We assume that
\begin{equation}
r(x) \in C^{\infty}(\dS^1), \quad r(x) > 0, \quad |r'(x)|\leq 1.
\end{equation}
The Laplace-Beltrami operator on $M$ is
\begin{equation}
\frac{1}{b^2}\left(-\frac{1}{r^m(x)}\partial_x r^m(x)\partial_x -
\frac{1}{r(x)^2}\Delta_Y\right).
\end{equation}
Letting $E_{\nu}$ be one of the eigenvalues of $-\Delta_Y$, we
consider  the operator $\D_{\nu,b}$ on $\dS^1$ given by
\begin{equation}
- \D_{\nu,b} = \frac{1}{b^2}\left(-\frac{1}{r^m(x)}\partial_x
r^m(x)\partial_x + \frac{E_{\nu}}{r^2(x)}\right).
\end{equation}
Let us note that one can compute the constant $b$ from the
asymptotics of eigenvalues of the operator $- \D_{\nu,b}$.
Therefore,  below we take $b=1$ for the sake of simplicity and let
$\D_{\nu}=\D_{\nu,1}$. Our problem is now stated as follows:

\begin{quote}
{\it Determine $r(x)$ from the knowledge of the spectrum of $-
\D_{\nu}$.}
\end{quote}

Assume that  $r(x)$ is 1-periodic and  is given by
\[
\lb{2}
 r(x)=r_0e^{{2\/m}Q(x)},\qq Q(x)=\int_0^xq(t)dt,\qq x\in [0,1],
\]
where $r_0 > 0$ is a constant. Our results are stated in terms of
$q$ and summarized in Section 2.

\subsection{Analytic approach to inverse problems}

There are various methods for inverse problems. We  employ here the analytic (or direct) approach due to  \cite{GT87},
\cite{KK97} based on nonlinear functional analysis, which we briefly explain below. Suppose that
$\cH,  \cH_1 $ are real separable Hilbert spaces. The derivative of
a map $f: \cH \to \cH_1 $ at a point $y\in \cH$ is a bounded linear
map from $\cH$ into $\cH_1$, which we denote by ${d\/dy}f$. A map
$f: \cH \to \cH_1 $ is compact on $\cH$, if it maps any weakly
convergent sequence in $\cH$ to a strongly convergent sequence in
$\cH_1.$  A map $f: \cH \to \cH_1 $ is a real analytic isomorphism
between $\cH$ and $\cH_1$, if $f$ is bijective
 and both $f$ and $f^{-1}$ are real analytic.
Let $\cH_C$ be the complexification of $\cH$.
 We recall the main result in \cite{KK97}, \cite{KK99}.

\begin{theorem}
\lb{TA97}
 Let $\cH, \cH_1$ be real separable Hilbert spaces equipped with norms
$\|\cdot \|, \|\cdot \|_1$. Suppose that a map $f: \cH \to \cH_1$
satisfies the following conditions:

\no (A) $f$ is real analytic,

\no (B)  ${d \/dq}f$ has an inverse for each fixed $q\in \cH$,

\no (C) there is a nondecreasing function $\e: [0,\iy ) \to [0,\iy
)$ such that $\eta(0) = 0$, and
$$
\|q\|\le \e (\|f(q)\|_1), \qqq   \forall \ q\in
\cH.
$$
\no (D) there exists a linear isomorphism $f_0:\cH\to \cH_1$ such
that the mapping $f-f_0: \cH \to \cH_1$  is compact.


 Then $f$ is a real analytic isomorphism between $\cH$ and $\cH_1$.
\end{theorem}

In our applications, $f(q)$ is supposed to be a spectral
data. The main issue is to derive a priori estimates of $\|q\|$ in terms of
spectral data $\e (\|f(q)\|_1)$.


\medskip

\section {Main results}
\setcounter{equation}{0}

\subsection{Inverse problem for perturbed torus.}

We assume that the profile $r(x), x\in [0,1]$, satisfies (\ref{2}) with
\[
\lb{f2}
q\in \mH_0\cap \mH_1,
 \quad
r(0)=r(1)=r_0>0,
\]
for a fixed radius $r_0$.  Here $\mH_j$ is the Sobolev space of real
functions defined by
$$
\mH_j=\mH_j(\T)=\rt\{q\in L^2(\T)\,  ; \, q^{(j)}\in L^2(\T),\, \int
_0^1q(x)dx=0\rt\},\qqq j\ge 0,
$$
equipped with  norm
$$
\|q\|^2_j=\|q^{(j)}\|^2=\int_0^1|q^{(j)}(x)|^2dx,
$$
and
$$
\mH_{-1}=\mH_{-1}(\T)=\rt\{f=g'\, ; \, g\in \mH_0\rt\}
$$
equipped with  norm
$$
\|f\|^2_{-1}=\|g\|^2=\int_0^1|g(x)|^2dx.
$$
We also
introduce a real Hilbert space $\ell^2_{j}$ of sequences $f=(f_n)_1^{\iy }$
equipped with  norm
\[
\lb{1.0} \|f\|_{j}^2=\sum _{n\ge 1}(2\pi n)^{2j}|f_n|^2,\qqq  j\in
\R.
\]

We  fix $\nu$ arbitrarily and omit the subscript $\nu$ of  $- \D_\n$
for the sake of simplicity. In \cite{K00}, \cite{KK00}, we have
already developed  the inverse spectral theory for the case
$E_\n=0$. Below we consider the case $E_\n>0$.   The spectrum of $-
\D_\n$ on $\T$ is discrete   and consists of eigenvalues
$\l_{2n}^\pm=\l_{\n,2n}^\pm, n\ge 0$, for the equation
\[
\lb{eq1} -{1\/\r^2}(\r^2 f')'+{E_\n\/r^2} f=\l f
\]
with periodic boundary condition $y(x+1) = y(x), x \in {\R}$.
They are  labeled as
\[
\l_0^+<\l_2^-\le \ \l_2^+\ <\l_4^-\le \ \l_4^+\ <\cdots.
\]
The anti-periodic eigenvalues
$\l_{2n-1}^\pm, n\ge 1$, are the eigenvalues  of the equation \er{eq1}
under antiperiodic boundary condition,  i.e., $y(x+1)=-y(x),\, \
x\in \R$. Thus the eigenfunctions corresponding to $\l_n^{\pm}$ are periodic with
period 1 for even $n$, and  antiperiodic for odd $n$.
It is well-known that the periodic eigenvalues $\l_0^+, \
\l_{2n}^\pm, n\ge 1$, determine the anti-periodic
eigenvalues $\l_{2n-1}^\pm, n\ge 1$, and that they satisfy
\[
\lb{psp1} \l_0^+<\l_1^-\le \ \l_1^+\ <\l_2^-\le \ \l_2^+<\l_3^-\le \
\l_3^+\ <\l_4^-\le \ \l_4^+\ <\cdots.
\]
The gap $\g_n$ for the operator $- \D_\n$ is defined by
\[
\lb{gap1} \g_n=[\l_n^-,\l_n^+],\qqq n\ge 1.
\]

We also use the spectral data for the Sturm-Liouville problem in the impedance
 form under Dirichlet boundary condition:
\[
\lb{Dbc} -{1\/\r^2}(\r^2 f')'+{E_\n\/r^2} f=\l f, \qqq f(0)=f(1)=0.
\]
Let  $\m_n=\m_n(q), n\ge 1$, and $f_n=f_n(x,q)$ be the eigenvalues
and the associated eigenfunctions. It is well-known that all $\m_n$
are simple and are labeled as $ \m_1<\m_2<\m_3< \cdots$. Moreover,
the Dirichlet eigenvalue $\m_n$ belongs  to the interval $\g_n$,
i.e.,
\[
\lb{dgap}
\m_n\in [\l_n^-,\l_n^+], \qqq \forall\ n\ge 1.
\]

As in \cite{K99}, \cite{K03}, we  construct  the
gap length mapping $\p$ by
\[
\lb{gm}
\p : q\to \p(q) =(\p_n)_1^\iy,\qqq \p_n=(\p_{n1},\p_{n2})\in \R^2,
\]
where $\p_{n1}$ and
$\p_{n2}$ are given by
\[
\begin{aligned}
\p_{n1}={\l_n^-+\l_n^+\/2}-\m_n,\qqq
\p_{n2}=\rt| {|\g_n|^2\/4}-\p_{n1}^2\rt|^{1\/2}\sign \vk_n,
\end{aligned}
\]
 the norming constants $\vk_n=\vk_{\n,n}$ are defined by
\[
\lb{f4} \vk_n=\log \rt|{f_n'(1)\/f_n'(0)}\rt|,\qqq n\ge 1,
\]
and where $f_n$ is an eigenfunction for the Dirichlet problem
\er{Dbc} associated with the eigenvalue $\m_n$. Such mapping was
introduced in \cite{K99} for the Schr\"odinger operator   with
periodic potential on the circle. It is known that $\p$ maps $\mH$
into $ \ell^2\os \ell^2$ (see \cite{K02}). Note that the vector $\p$
depends on  the gap lengths $|\g_n|$, ${\rm sign }\, \vk_n$ and
$\p_{n1}$ for all $n\geq 1$, but is independent of the location of
the gaps and the Dirichlet eigenvalues themselves.

The standard arguments (see \er{5}-\er{6}) show that the operator
$\D_\n$ is unitarily equivalent to a one-dimensional Schr\"odinger
operators $S_\n y=-y''+p y$ on $L^2(\T)$, where the potential $p$ is
given by
\[
\begin{aligned}
\lb{62}  p=q'+q^2+u_\n-c_{\n,0},\qq u_\n={E_\n\/r_0^2}e^{-{4\/m}Q},
 \qq
 c_{\n,0}=\int_0^1(q^2+u_\n)dx.
\end{aligned}
\]
If $q\in \mH_j, j\ge 1$, then we deduce that the potential $p\in
\mH_{j-1}\ss L^2(\T)$. If $q\in \mH_0$, then we deduce that the
potential $p\in \mH_{-1}$ is a distribution.

\begin{theorem}
\lb{T1} Fix $\n\ge 0$ and let $m\ge 1$. Then the mapping
$\p:\mH_1\to\ell^2\os\ell^2$, given by \er{gm}, is a  real
analytic isomorphism between the Hilbert spaces $\mH_1$ and
$\ell^2\os\ell^2$. Moreover, for all $q\in \mH_1$ the following
estimates hold true:
\[
\lb{13} \|q'\| \le\|p\| \le 2\|\p\|(1+\|\p\|^{1\/3}),
\]
\[
\lb{14} \|\p\|\le w(1+w^{1\/3}),\ \ \
\]
where $w=\|q'\|+\|q\|(\|q'\|+ C_*^{1\/2} e^{\b\|q\| })$ and
$$
\begin{aligned}
 C_*=A(\b+\a)(2+\b A),\qqq \b={4\/m}, \qqq A={E_\n\/r^2_0}.
\end{aligned}
$$
\end{theorem}

We consider the inverse  problem for more singular  coefficients
$q\in \mH_0$. In this case, in general, the gap length $|\g_n|$ is
increasing as $n\to \iy$.

\begin{theorem}
\lb{T2} i) Fix $\n\ge 0, m\ge 1$. Then the mapping
$\p:\mH_0\to\ell_{-1}^2\os\ell_{-1}^2$, given by \er{gm}, is real
analytic, injective  and  satisfies
\[
\lb{13x}
\begin{aligned}
\|\p\|_{-1}\le \|p\|(1+2\|p\|)^3,
\\
\|p\|\le   \|q\|(3+ 2\|q\|+\b Ae^{\b \|q\|}).
\end{aligned}
\]
ii) Let, in addition, $m=1$. Then the mapping
$\p:\mH_0\to\ell_{-1}^2\os\ell_{-1}^2$ is a real analytic
isomorphism between the Hilbert spaces $\mH_0$ and
$\ell_{-1}^2\os\ell_{-1}^2$ and
\[
\lb{14x}
\begin{aligned}
 \|q\|^2\le 2\|p\|^2(1+2\|p\|^2),
\\
\|p\|\le 96\pi^2 \|\p\|_{-1}(1+2\|\p\|_{-1})^3.
\end{aligned}
\]

\end{theorem}
{\bf Remark.} 1) If we know the location of  gaps $\g_n$, Dirichlet
eigenvalues $\m_n$ and  $\sign \vk_n$ for all  $n\ge 1$, then we can
obtain $q$ uniquely. Moreover,  in order to determine $q$ we need to
know only the vector $\p(q)$. From the value $\p(q)$, we can compute
the gap length $|\g_n|$, the distance $\p_{cn}$ between the
Dirichlet eigenvalue and the center of the gap for any $n\ge 1$.

2) The proof of Theorems \ref{T1} and \ref{T2}  is based on
non-linear functional analysis. It is crucial that the mapping
$\p(q)$ is a composition of two mappings $\P(p)$ and $p=P(q)$ given
by \er{gmx} and \er{dP1}. All necessary properties for $\P(q)$ are
described in Theorem \ref{Tgap}. The mapping $q\to P(q)$ is defined
by Theorem \ref{T3}. The a priori estimates of $\|q\|$ in terms of
$\|P(q)\|$ play an  important role in the  proof.

\subsection {Perturbed Riccati mappings}
Recall that $q\to q'+\a q^2$ is the Riccati mapping. In order to
prove our main theorems we use perturbed Riccati mappings
$P:\mH_j\to \mH_{j-1}$  given  by
\[
\begin{aligned}
\lb{dP1} & q\to P, \qqq
 P(q)=q'+\a q^2+u-c_0,\\
& u=Ar^{-\b Q},\qqq
\qq Q(x)=\int_0^xq(t)dt,\qq
 c_0=\int_0^1(\a q^2+u)dx,
\end{aligned}
\]
where $\a, \b>0$ and $A\ge 0$. We also consider the singular case
$q\in \mH_0$, which is  more complicated.

\begin{theorem}
\lb{T3}  i) The mapping $P:\mH_j\to \mH_{j-1}, j\ge 1$, given by
\er{dP1} is a real  analytic isomorphism between the Hilbert spaces
$\mH_j$ and $\mH_{j-1}$ and satisfies:
\[
\lb{t3e} \|q'\|^2\le  \|P(q)\|^2 \le \|q'\|^2+\a ^2
2\|q\|^3\|q'\|+C_*\|q\|^2e^{2\b \|q\|},\qqq
\]
where the constant $C_*=A(\b+\a)(2+\b A)$.

\noindent
ii) Let $j=0$. Then the mapping $P:\mH_0\to \mH_{-1}$, given by
\er{dP1} is real  analytic and  injective. Moreover,
\[
\lb{pq01} \|P(q)\|_{-1}\le   \|q\|(3+ 2\|q\|+\b Ae^{\b \|q\|}).
\]
If $m=1$, then the mapping $P:\mH_0\to \mH_{-1}$ is a real analytic
isomorphism between the Hilbert spaces $\mH_0$ and $\mH_{-1}$ and
satisfies:
\[
\lb{t4e}  \|q\|^2\le 2\|P(q)\|_{-1}^2(1+2\|P(q)\|_{-1}^2) .
\]
\end{theorem}

\noindent {\bf Remark.} 1) When $m \geq 2$,  we do not have
estimates of $\|q\|$ in terms of $\|P(q)\|$, which causes a
difficulty in the inverse problem for the mapping $P:\mH_0\to
\mH_{-1}$.

\subsection {Symmetric surfaces}

Define the spaces of even functions $\mH_j^{even}(\T)$,  and of odd
functions $\mH_j^{odd}(\T)$  by
\[
\begin{aligned}
\lb{oeL} \mH_j^{even}&=\rt\{q\in \mH_j: \, q(x)=q(1-x), \qq \forall
\ x\in (0,1)\rt\},\\
\mH_j^{odd}&=\rt\{q\in \mH_j: q(x)=-q(1-x), \qq \forall \ x\in
(0,1)\rt\}.
\end{aligned}
\]
Note that we have $\mH_j=\mH_j^{even} \os \mH_j^{odd}$.
It is well-known that if $P(q)\in \mH_j^{even}$,  we have either
$\m_n=\l_n^- $ or $\m_n=\l_n^+ $, see \cite{GT84}. Then, $\p_{n2}=0$ and
$|\p_{n1}|=|\g_n|/2$ for all $n\ge 1$.

\begin{corollary}   \lb{Todd}
i) The mapping $P:\mH_j^{odd}\to \mH_{j-1}^{even}, j\ge1$,  given by
\er{dP1} is a real analytic isomorphism between the Hilbert spaces
$\mH_j^{odd}$ and $\mH_{j-1}^{even}$.

\noindent
ii) Fix $\n\ge 0$ and let $m\ge 1$. Then the mapping
$\p^e:\mH_1^{odd}\to\ell^2$, given by
\[
\begin{aligned}
\lb{epj} q\to \p^e=(\p_{n}^e)_1^\iy, \qq
\p_{n}^e={\l_n^-+\l_n^+\/2}-\m_n
\end{aligned}
\]
 is a  real analytic isomorphism between the Hilbert spaces
$\mH_1^{odd}$ and $\ell^2$.

Moreover, if $m=1$ then the mapping
$\p^e:\mH_0^{odd}\to\ell_{-1}^2$, given by \er{epj}
 is a  real analytic isomorphism between the Hilbert spaces
$\mH_0^{odd}$ and $\ell_{-1}^2$.
\end{corollary}

\noindent {\bf Remark.}  Similarly to \er{13} and \er{14}, one can
also derive the estimates of $q$ in terms of gap-length.

\subsection {Minkowski problem }

The Minkowski problem deals with the
existence of a convex surface with a prescribed Gaussian curvature.
More precisely, for a given strictly positive real function $F$
defined on a sphere, one seeks a strictly convex compact surface
$\cS$, whose Gaussian curvature at $x$ is equal to $F({\bf n}(x))$,
where ${\bf n}(x)$ denotes the outer unit normal to $\cS$ at $x$.
The Minkowski problem was solved by Pogorelov \cite{P74} and by
Cheng-Yau \cite{CY76}.

We consider only the case $m=\dim Y=1$. Note that our surface is not
convex, in general.
 We solve an analogue of the Minkowski problem
in the case of the surface of revolution by showing the existence of
a bijection between the  Gaussian curvatures and the profiles of
surfaces.

It is well-known that the Gaussian  curvature $G$ is given by
\[
\lb{K1} G=-{r''\/r}=-v'-v^2,\qqq v=2q,\qq \r=r^{1\/2},\qq {\rm for} \quad
m=1.
\]
We define a new variable $G_1\in \mH_0$ by
\[
\lb{b2}
\begin{aligned}
G=G_0+G_{1}, \qqq G_{0}=\int_0^1 G(x)dx,
\\
G_1=-v'-v^2+\int_0^1v^2dx.
 \end{aligned}
\]
Assuming that the function $G_1\in \mH_0$ is given, we determine
the profile (radius) $r=r(x), x\in
[0,1]$, and the constant $G_0$ by solving the
equation \er{K1}.


Recall that $m=\dim Y$. The eigenvalues $\cE$ and  $e_j,  j=1,
\cdots, m$, of the Ricci tensor $({\rm
Ric}_{ij})_{i,j=0,1,..\cdots,m}$ are given by
\[
\lb{REi}
\begin{aligned}
& \cE= -v'-{v^2\/m},\qqq v=2q,
 \\
&  e_j=
  \frac{\kappa_j}{r^2} - {1\/m}(v'+v^2),\qq \qqq j=1,...,m,
  \end{aligned}
\]
 where $\kappa_j>0$ denote the eigenvalues of the
Ricci tensor $({\rm Ric}_{ij}^Y)_{i,j=1,\dots,m}$ on $Y$.  In
particular, if $Y$ is a sphere, then all eigenvalues
$\kappa_j=\kappa>0$, $j=1,\dots m$ for some $\kappa>0$.  In this
case, the Ricci tensor of the warped product has one simple
eigenvalue $E$ and an eigenvalue $e_1=...=e_m$ of multiplicity $m$
given by
\[
\lb{RiSp} \cE=-v'-{1\/m}v^2,\qqq  \qq
e_1={\kappa\/r^2}-{1\/m}(v'+v^2).
\]
Note that if $m=1$, then the eigenvalues $\cE$ of the Ricci tensor
coincides with the Gaussian  curvature $G$, i.e.,
\[
\cE=G\qqq  {\rm at} \qq m=1.
\]
For the function $q\in \mH_j$ we define the constant $\cE_0$ and the
function $\cE_1\in \mH_{j-1}$ by
\[
\lb{Rcr1}
\begin{aligned}
& \cE=\cE_0+\cE_1,\qqq \cE_{0}=\int_0^1\cE dx
=-{1\/m}\int_0^1v^2(x)dx\le 0,
 \\
&   \cE_1=-v'-{1\/m}v^2+{1\/m}\int_0^1v^2(x)dx.
\end{aligned}
\]
If $v=2q\in \mH_j$, then \er{Rcr1} gives that $\cE\in \mH_{j-1}$.
The function $\cE_1(x), x\in [0,1]$, is the non-constant part of
eigenvalue $\cE(x), x\in [0,1]$ of the Ricci curvature tensor and
$E_{0}=\const$ is the constant part of eigenvalue $\cE$.


\begin{corollary}
\lb{T4} Let $m\ge 1$.  Then the mapping $q\to \cE_1$ given by
\er{Rcr1} is a real analytic isomorphism between the Hilbert spaces
$\mH_1$ and $\mH_0$ and satisfies:
\[
\lb{13c} \|v'\|^2\le \|\cE_1\|^2\le
\|v'\|^2+{\|v\|^2\|v'\|^2\/m^2}-{\|v\|^4\/m^2}.
\]
 Moreover, the constant $\cE_0$ is uniquely defined by $\cE_1$.
\end{corollary}

\noindent {\bf Remark.} 1) In the case $m=1$, this corollary solves
the Minkowski problem for the perturbed torus.

\noindent 2) This theorem gives an isomorphism between $\cE_1$ and
the function $q$. Thus, in addition to Theorem \ref{T1}, we can get
a parametrization of the surface by $\cE_1$ or  the function $q$ or
the vector $\p$:
$$
q\qq \Longleftrightarrow  \qq \cE_1  \qq \Longleftrightarrow  \qq
\p.
$$


\subsection{Brief overview}
There is an abundance of works devoted to the spectral theory and
inverse problems for the surface of revolution from the view points
of classical inverse Strum-Liouville theory, integrable systems,
micro-local analysis, see \cite{AA07}, \cite{E98}
and references therein. For integrable
systems associated with surfaces of revolution, see e.g.
\cite{KT96}, \cite{Ta97},
\cite{SW03}
and references therein.

Bruning-Heintz \cite{BH84} proved that the symmetric metric is
determined from the spectrum by using the standard 1-dimensional
Gel'fand-Levitan theory. We mention the work of Zelditch \cite{Z98},
which proved that the isospectral revolutionary surfaces of simple
length spectrum, with some additional conditions, are isometric.
Isozaki-Korotyaev \cite{IK17}, \cite{IK17x}  solve the inverse
spectral problem for rotationally symmetric manifolds (finite
perturbed cylinders), which include a class of surfaces of
revolution, by giving an analytic isomorphism from the space of
spectral data onto the space of functions describing the radius of
rotation. An analogue of the Minkowski problem is also solved.

As far as we know, there were no results about inverse problems for perturbed torus.

 In this paper we use inverse spectral theory for  Schr\"odinger operators  with
potentials on the circle. Let us briefly review the inverse spectral
theory for Strum-Liouville operators on the circle, mostly focusing on the {\it characterization}
problem, i.e., the complete description of spectral data that
correspond to some fixed class of potentials. More information about
different approaches to inverse spectral
 problems can be found in the monographs \cite{M86}, \cite{L87} and the papers
 \cite{MO75}, \cite{K97}, \cite{K99} and \cite{K03}
 and references therein.

Dubrovin \cite{D75},
Its and Matveev \cite{IM75},  Novikov \cite{N74}, Trubowitz
\cite{T77}  considered the inverse problem for finite band
potentials (potentials were more general in \cite{T77}).
Marchenko-Ostrovski \cite{MO75} solved the inverse problem including
the characterization in terms of spectral data associated with a
"global quasimomentum". Note that their construction and the proof
are complicated, for example, they used the inverse spectral theory
for the scattering on the half-line.

Garnett--Trubowitz \cite{GT87} solved the inverse problem
for the gap-length mapping for the case of even potentials, under
the conjecture that there exists an estimate of  potentials in terms
their gap lengths.
 Kargaev--Korotyaev \cite{KK97} gave a simplified proof of  the result of Garnet and Trubowitz
\cite{GT87}. The Garnett--Trubowitz conjecture on the estimate of
potentials in terms of their gap lengths was proved by Korotyaev
\cite{K98}, \cite{K00}. Korotyaev \cite{K97} reproved shortly
Marchenko-Ostrovski results \cite{MO75}. Korotyaev \cite{K99}
constructed the gap-length mapping given by \er{gmx} for potentials
from $\mH_0$ and from $\mH_{-1}$  and solved the corresponding
inverse problems. Note that for the case $\mH_{-1}$ inverse problems
for the operator in impedance form are important and  were solved
for the Dirichlet problem Coleman-McLaughlin \cite{CM93} and for
periodic case by Korotyaev \cite{K99}.

We use also results on perturbed Riccati mappings from
\cite{K02}, \cite{K03}, \cite{BKK03}.

\subsection{Plan of the paper}
In Section 2, we prove Theorems \ref{T1} and \ref{T3}, which are based on the
 Theorem \ref{TA97} in non-linear functional analysis
\cite{KK97}. We do it after preparing a priori estimates for the
perturbed Riccati mapping. In Section 3 we prove Theorems \ref{T2},
where the main problem is a priori estimates for the perturbed
Riccati mapping. For the moment we can do it only for  $m=1$. In
Section 3, we consider the mapping associated with the Minkowski
problem  and the eigenvalue of the Ricci tensor.


\section {Proof of Theorem \ref{T1}}
\setcounter{equation}{0}

\medskip

\subsection {The  unitary transformations.}
We begin by explaining the crucial role of the factorization  and the
non-linear mapping $P$ given by \er{dP1}. For each $q\in \mH_j, j\ge 0$, we define the periodic
weighted (model) operator  $-\D_\n $ given by
\[
\lb{1.9u}
 -\D_\n f=-{1\/r^m}(r^m f')'+u_\n f, \qq u_\n={E_\n\/r^2},
 \qq  r=r_0e^{{2\/m}Q},\qq Q= \int_0^xq(t)dt,
\]
where $f\in L^2(\T,r^mdx)$. The property of the operator $-\D_\n $ is well-known
(see \cite{Kr57}, \cite{K00}, \cite{K03}). The
spectrum of the operator $-\D_\n $ is discrete and given by
\er{psp1}. We define the unitary transformation $\mU:
L^2(\T,\r^2dx)\to L^2(\T,dx)$ by
 $$
\mU f= \r f,\qq \r=r^{{m\/2}}=\r_0e^Q,\qq \r_0=r_0^{{m\/2}}.
$$
Using $\mU$ and the equality
$r=r_0e^{{2\/m}\int_0^xq(t)dt}$, we transform the operator $- \D_\n$ into the
Schr\"odinger operator $S_\n$:
\[
\lb{5}
\mU (- \D_\n) \mU^{-1}=
 -\r^{-1}\pa_x \r^2\pa_x \r^{-1}+u_\n=\cD^*\cD+u_\n=S_\n+c_{\n,0},
\]
where
\[
\begin{aligned}
\lb{a7} \cD=\r\ \pa_x \ \r^{-1}= \pa_x -q,\qqq \cD^*=-\pa_x-q,\\
\cD^*\cD=-(\pa_x +q)(\pa_x -q)=-\pa_x^2+q'+q^2,
\end{aligned}
\]
\[
\lb{a7x} S_\n y =-y''+P_\n(q)y,
\]
acting on $ L^2(\T,dx)$, and
where the periodic potential $P_\n(q)$ has the form \er{dP1} with
$A={E_\n\/r_0^2}, \b={4\/m}, \a=1$, i.e.,
\[
\begin{aligned}
\lb{6}  P_\n(q)=q'+q^2+u_\n-c_{\n,0},\qqq
u_\n={E_\n\/r_0^2}e^{-{4\/m}Q},
 \qqq
 c_{\n,0}=\int_0^1(q^2+u_\n)dx.
\end{aligned}
\]
Thus we see that the Laplacian on $(M, g_M)$ is unitarily equivalent
to a direct sum of one-dimensional Schr\"odinger operators $S_\n$,
namely,
\[
\lb{} \D_{(M,g)}\backsimeq \os_{\n\ge 1} (S_\n+c_{\n,0}),
\]
where the direct sum acts on $\os_{\n\ge 1} L^2(\T,dx)$.

\subsection{Schr\"odinger operators  on the circle}
  Using the global transformation $q\to P_\n(q)$ we reduce our operators $-\D_\n$ to the
Schr\"odinger operators $S_\n=-{d^2\/dx^2}+ P_\n(q)$ with the
potential $P_\n(q)$ (including the singular case $P_\n(q)\in
\mH_{-1}$).

We recall some results about inverse spectral theory for
Schr\"odinger operator with potential $p\in \mH_{j}, j\ge -1$. We
consider the operator $S= -{d^2\/dx^2}+p(x),$ acting on $L^2(\T)$,
 where the 1-periodic  potential $p$ belongs to the Hilbert space
$\mH_j, j \ge -1$.   The spectrum of $S$ on $\T$ is discrete
  and consists of eigenvalues  $\l_{0}^+$ an $\l_{2n}^\pm, n\ge 1$,
  which are eigenvalues of the equation
\[
\lb{se} -f''+p f=\l f.
\]
with 1-periodic boundary conditions. They are  labeled by
\[
\lb{A2z} \l_0^+<\l_2^-\le \ \l_2^+\ <\l_4^-\le \ \l_4^+\ <\dots
\]
It is well known that if we know all periodic eigenvalues $\l_0^+, \
\l_{2n}^\pm, n\ge 1$, then we can recover so-called anti-periodic
eigenvalues $\l_{2n-1}^\pm, n\ge 1$. The anti-periodic eigenvalues
$\l_{2n-1}^\pm, n\ge 1$ are eigenvalues of the equation \er{se}
under antiperiodic boundary conditions,    i.e., $y(x+1)=-y(x),\ \
x\in \R$. Thus the eigenfunctions corresponding to $\l_n^{\pm}$ have
period 1 when $n$ is even and they are antiperiodic, when $n$ is
odd. It is important that they satisfy
\[
\l_0^+<\l_1^-\le \ \l_1^+\ <\l_2^-\le \ \l_2^+<\l_3^-\le \ \l_3^+\
<\l_4^-\le \ \l_4^+\ <\dots
\]
Introduce so-called gaps $\g_n$ for the operator $S$ by
\[
\g_n=[\l_n^-,\l_n^+],\qqq n\ge 1.
\]
It is important that each Dirichlet eigenvalue $\m_n, n\ge 1$
belongs to the interval $[\l_n^-,\l_n^+]$, i.e.,
\[
\lb{dgapx} \m_n\in [\l_n^-,\l_n^+] \qqq \forall\ n\ge 1.
\]
Introduce the fundamental solutions $\vp(x,\l,p)$ and $\vt(x,\l,p)$
of the equation
\[
\lb{eqp} -y''+py=\l y, \ \ \ \l\in \C,
\]
under conditions
  $$
  \vp'(0,\l,p)= \vt (0,\l,p)=1,\qqq \qqq
\vp(0,\l,p)=\vt'(0,\l,p)=0.
$$
We introduce the Lyapunov function (the discriminant) $\L$ by
$$
\L(\l,p)={1\/2}(\vp '(1,\l,p)+\vt(1,\l,p)).
$$
 Note that
\[
\begin{aligned}
&\L(\l_{n}^{\pm},p)=(-1)^n,\qqq \forall \  n\geq 1. \\
\end{aligned}
\]
Here and in the sequel, $  ' =\partial /\partial x , \ \
\dot {\  }=\partial /\partial z, \ \ \ \partial =\partial /\partial q.$

\subsection{Gap length mapping}
Using the gaps $\g_n=[\l_n^-,\l_n^+], n\ge 1$, the Dirichlet
spectrum $\m_n(p)$ and the sign of the  norming constants $\sign
K_n$, we now construct the gap length mapping $\P: \mH_1\to
\ell^2\os \ell^2$ by
\[
\lb{gmx} p\to \P=(\P_n)_1^\iy,\qqq \P_n=(\P_{n1},\P_{n2})\in \R^2.
\]
Here the coordinates $\P_{n1},\P_{n2}$ are given by
\[
\lb{P12}
\begin{aligned}
\P_{n1}={\l_n^-+\l_n^+\/2}-\m_n, \qqq \P_{n2}=\rt|
{|\g_n|^2\/4}-\P_{n1}^2\rt|^{1\/w}\sign K_n,
\end{aligned}
\]
where the norming constant $K_n=\log |\varphi '(1,\m_n,p)|$ is
defined by \er{f4}. Such mapping was introduced in \cite{K99} for
the Schr\"odinger operator with periodic potential. It is known that
$\P$ maps $\mH_j$ into $ \ell_{j}^2\os \ell^2$ (see \cite{K99}).
Note that  $\P$ is computed from  the gaps lengths
$|\g_n|, {\rm sign } K_n$ and $\P_{n1}$ for any $n\ge 1$, however,
we do not need the position of the gaps and the Dirichlet
eigenvalues.

 We introduce the
Fourier transformation $\Phi : H_C\to \ell^2_C$ by $(\Phi
q)_n=\sqrt{2} q_n, n\geq 1, $ where $q_n=(q_{cn}, q_{sn})$ with
$$
 q_{cn}=\int _0^1 q(x)\cos 2\pi nx dx,\ \
 q_{sn}=\int _0^1 q(x)\sin 2\pi nx dx,\ \  n\geq 1,\ \ q\in H_C.
$$

We now formulate the needed results from \cite{K99}, \cite{K98} for
$j=0$ and from \cite{K03} for $j=-1$, where the inverse problem for
the mapping $ \P(\cdot )$ is solved by the direct method.

\begin{theorem}
\lb{Tgap} Each mapping $\P: \mH_j \to \ell_{j}^2\oplus \ell_{j}^2,
j=-1, 0$ defined by \er{gmx}, is a real analytic isomorphism between
$\mH_j$ and $\ell_{j}^2\oplus \ell_{j}^2$.
 Moreover the following estimates hold true:
\[
\begin{aligned}
& \|p\| \leq 2\|\P\|(1+\|\P\|^{1\/3}), \\
 & \|\P\| \leq\|p\|(1+\|p\|^{1\/3}),
\end{aligned} \qqq {\rm if}\qq j=0,
\]
and
\[
\lb{1.7}
\begin{aligned}
&  \|\vk\|\le 96\pi^2 \|\P\|_{-1}(1+2\|\P\|_{-1})^3,\\
&  \|\P\|_{-1}\le \|\vk\|(1+2\|\vk\|)^3.
\end{aligned}\qqq {\rm if}\qq j=-1,\qq p=\vk', \vk\in \mH_0.
\]
 \end{theorem}

\no  {\bf Remark.}  The proof is based on on nonlinear functional
analysis see Theorem \ref{TA97}.  The Gelfand-Levitan-Marchenko
equation and a trace formula are not used in the proof.

\subsection{Proof of Theorem \ref{T1}}
We assume that Theorem \ref{T3} i) holds true (in our case $\a=1$)
and then the mapping $q\to p=P_\n(q)$ given by \er{6} is a
real-analytic isomorphism between $\mH_1$ and $\mH_0$. Recall that
due to  Theorem \ref{Tgap} from \cite{K99} the mapping
\[
\lb{FKC} \P: p\mapsto \P(p)= (\P_n(p))_{n=1}^{\iy}
\]
is a real-analytic isomorphism between $\mH_1$ and
$\ell_2\os\ell_2$.

Thus due to \er{gm}--\er{f4}, Theorem \ref{Tgap} and the relation
between $q$ and $P$ given by \er{a7x}, \er{6}, we obtain the identity
$$
\p(q)=\P(P(q)),\qq \forall \ q\in \mH_1.
$$
The mapping $\p(\cdot)$ is the composition of two mappings $\P$ and
$P$, where each of them is the corresponding analytic isomorphism
(see Theorem \ref{Tgap} and Theorem \ref{T3}). Then the mapping
$$
\p:q\mapsto (\P_n(P(q)))_{n=1}^{\iy},
$$
is a real-analytic isomorphism between $\mH_1$ and
$\ell^2\os\ell^2$. Using the factorization  $\p(q)=\P(P(q)) $ and
combining the estimates in Theorems \ref{T3} and \ref{Tgap}, we
obtain
$$
\|P\| \le 2\|\p\|(1+\|\p\|^{1\/3}), \qqq \|\P\|
\leq\|P\|(1+\|P\|^{1\/3}),
$$
$$
\|q'\| \le \|P\|\le \|q'\|+ 2\|q\|\|q'\|+C_*^{1\/2}\|q\|e^{\b \|q\|},
$$
where $C_*=A(\b+\a)(2+\b A)$ and $\b={4\/m}, A={E_\n\/r^2_0}$. This
yields \er{13}, \er{14}. \BBox

\section {Proof of Theorem \ref{T3} for smooth case }



 \subsection {A priori estimates}
 Define the scalar product in $L^2(0,1)$ by
$$
(f,g)=\int_0^1f\ol gdx.
$$
We recall the following simple estimates
\[
\lb{qdq}
 \|q\|_{L^\iy(0,1)} \le  \|q'\|, \qqq \|q\| \le  \|q'\|, \qqq
 \forall \ q\in \mH_1,
\]
and
\[
\begin{aligned}
\lb{Qq}
\|Q\|_{L^\iy(0,1)} \le  \|q\|,\quad
 \|u(Q)\|_{L^\iy(0,1)} \le Ae^{\b \|q\|}.
 \end{aligned}
\]

We consider the mapping $q\to P(q)$ given by
\[
\begin{aligned}
\lb{dP}
& P(q)=q'+\a q^2+u-c_0,\\
& \a>0,\qq   u(Q)=Ae^{-\b Q},\qq Q(x)=\int_0^xq(t)dt,\qq
 c_0=\int_0^1fdx,
\end{aligned}
\]
where  $f=\a q^2+u$ and   $A={E_\n\/r_0^2}$  and $\b={4\/m}$. We
derive two-sided estimates of $q$ and $P(q)$.

\begin{lemma}
\lb{Tpq1} Let the mapping $q\to P(q)$ from $\mH_j$ into $\mH_{j-1},
j\ge 1$, be given by \er{dP}. Then the following estimates hold true:
\begin{equation}
 \|q'\|^2\le  \|P(q)\|^2=\|q'\|^2+\| \a q^2+u-c_0\|^2+2\b
(q^2,u), \lb{Pe1}
\end{equation}
\begin{equation}
  \|P(q)\|^2=\|q'\|^2+\a ^2\|q^2\|^2+\|u\|^2+2(\b+\a)(q^2,u)
-c_0^2, \lb{Pe2}
\end{equation}
\begin{equation}
 \|P(q)\|^2 \le \|q'\|^2+\a ^2
2\|q\|^3\|q'\|+C_*\|q\|^2e^{2\b \|q\|}, \lb{Pe3}
\end{equation}
where the constant $C_*=A(\b+\a)(2+\b A)$.

In particular, if $u=0, \a=1$, i.e., $A=0$, then $c_0=\|q\|^2$ and
\[
\lb{Pe4} \|q'\|^2\le \|P(q)\|^2=\|q'\|^2+\|q^2-c_0\|^2
=\|q'\|^2+\|q^2\|^2-c_0^2 \le  \|q'\|^2+\|q^2\|^2.
\]
\end{lemma}
\no {\bf Proof.} Let $h=\a q^2+u-c_0$. We have $p'=q'+h, h=f-c_0$
and
$$
\begin{aligned}
&\|P\|^2=\|q'\|^2+\|h\|^2+2(q',h),\\
&(q',h)=(q',\a q^2+u)=(q',u)=\int_0^1q'u(Q)dx=\b \int_0^1q^2u(Q)dx=
\b (q^2,u),\\
 \end{aligned}
$$
where the integration by parts has been used. This yields \er{Pe1}.

By a direct computation
$$
\begin{aligned}
&\|h\|^2=\|\a  q^2+u-c_0\|^2=\|\a q^2+u\|^2-2(\a q^2+u,c_0)+c_0^2=\|\a q^2+u\|^2-c_0^2,\\
& \|\a q^2+u\|^2=\a ^2\|q^2\|^2+\|u\|^2+2\a (q^2,u),
\end{aligned}
$$
which, together with \er{Pe1}, yields \er{Pe2}.

We show \er{Pe3}. We have $c_0=\a \|q\|^2+u_0$ and $u=u_0+u_1$,
where $u_0=(u,1)$. Then
$$
\|u\|^2-c_0^2=\|u_1\|^2-\a ^2\|q\|^4-2\a \|q\|^2u_0\le \|u_1\|^2.
$$
We need to estimate $\|u_1\|^2$. Letting $Q_x=Q(x)$, we obtain
\begin{equation}
\begin{aligned}
\|u_1\|^2=\int_0^1 (u(Q_x)-u_0)^2dx,\quad
u(Q_x)-u_0=\int_0^1 (u(Q_x)-u(Q_t))dt,
\end{aligned}
\nonumber
\end{equation}
and for some $y\in [x,t]$
\begin{equation}
u(Q_x)-u(Q_t)=A\big(e^{-\b Q_x} -e^{-\b Q_t}\big)=-\b Ae^{-\b
Q_y}(Q_x-Q_t).
\nonumber
\end{equation}
These identities and the following estimates
$$
\begin{aligned}
|u(Q_x)-u(Q_t)|=A|e^{-\b Q_x} -e^{-\b Q_t}|\le \b \|q\| Ae^{\b
\|q\|}
\end{aligned}
$$
imply
$$
\|u\|^2-c_0^2\le \|u_1\|^2 \le \int_0^1 \b^2 \|q\|^2 A^2e^{2\b
\|q\|}dx=\b^2 \|q\|^2 A^2e^{2\b \|q\|}.
$$
Moreover, we have
$$
\|P\|^2=\|q'\|^2+\a ^2\|q^2\|^2+V,
$$
where,
$$
V=2(\b+\a )(q^2,u)+\|u\|^2 -c_0^2,
$$
which is estimated as follows:
\begin{equation}
\begin{aligned}
& V\le 2(\b+\a )\|q\|^2A e^{\b \|q\|}+\b^2 \|q\|^2A^2e^{2\b \|q\|}\\
& = A\rt(2\b+2\a+\b^2 A\rt)\|q\|^2e^{2\b \|q\|} \le A(\b+\a)(2+\b
A)\|q\|^2e^{2\b \|q\|}
\\
& = C_*\|q\|^2e^{2\b \|q\|},\qqq C_*=A(\b+\a)(2+\b A).
\end{aligned}
\nonumber
\end{equation}
In addition, we have
$$
\begin{aligned}
\|q^2\|^2=\int_0^1q^4(y)dy=2\int_0^1q^2(y)dy\int_0^y q(t)q'(t)dt\le
2\|q\|^3\|q'\|,
\end{aligned}
$$
since $q(x_*)=0$ for some $x_*\in [0,1]$ and we have used the new
variable $x=x_*+y$. Combining all these estimates, we obtain \er{Pe3}. If
$u=0$, \er{Pe4} follows from \er{Pe1}-\er{Pe3}. \BBox

\subsection { Analyticity and invertibility.} We show that that mapping
$P: \mH_j\to \mH_j, j\ge 0$, given by \er{dP},  is real analytic and
locally invertible. Here we have  $P(q)=P(q,A,\b,\a)$ for some fixed
$A,\b, \a$. First we discuss analyticity of  $P$. Let
$$
Jf=\int_0^xfdx.
$$

\begin{lemma}\lb{TAn}
i)  The map $P: \mH_j\to \mH_{j-1}, j\ge 0$, given by \er{dP}  is
a real analytic  and its gradient is given by
\[
\lb{A1}
 {\pa P(q)\/\pa q} f=f'+2q f -\b u(Q) J f-{\pa c_0(q)\/\pa q}f,
 \qqq \ \forall \ \ q,f\in \mH_1,
\]
\[
\lb{A2} {\pa c_0(q)\/\pa q}f=\int_0^1\rt[2qf -\b u(Q) J  f \rt]dx.
\]
\no ii) Moreover, for each $q\in\mH_j$ the linear operator ${\pa
P(q)\/\pa q} $ acting in $\mH_j$ is invertible.
\end{lemma}

\no {\bf Proof}. The standard arguments (see \cite{PT87}) give the proof of
 i).

ii)  Due to \er{A1} the linear operator ${\pa P(q)\/\pa q}: \mH_j\to
\mH_j$ is a sum of the linear operator operator ${d\/dx}$ and the
compact operator for all $p\in \mH_1$. Thus ${\pa P(q)\/\pa q}$  is
a Fredholm operator. We prove that the operator ${\pa P(q)\/\pa q}$
is invertible by contradiction. Let $f\in \mH_j$ be a solution of
the equation
\[
\lb{Cont} {\pa P(q)\/\pa q}f=0
\]
for some fixed $q\in \mH_j$. Setting $y=J f$, we obtain the following equation
for  $y$:
\[
\lb{EE1}
\begin{aligned}
-y''-2\a qy'+V y=B,\qqq B=\int_0^1 (-2\a qy' +V y)dx,\\
y\in \mH_{j+1},\qqq V=\b u(Q),\qqq y(0)=y(1)=0.
\end{aligned}
\]
We rewrite the equation \er{EE1} in the following form
\[
\lb{EE2} -{1\/\r^2}(\r^2y')'+V y=B,\qqq \a q={\r'\/\r},\qqq
y(0)=y(1)=0.
\]
We have 2 cases. Firstly, let $B=0$. Then  multiplying \er{EE2} by $\r^2y$ we have
$$
0=\int_0^1\rt[-(\r^2y')'y +\r^2 V y^2\rt]dx= \int_0^1\rt[(\r y')^2 +
\r^2V y^2\rt]dx,
$$
which shows $y = 0$.

Secondly, let $B\ne 0$. It is sufficiently to consider the case
$B=1$. The solution of the equation \er{EE2} with $B=1$ has the form
\[
\begin{aligned}
\lb{eB1} y=y_0+y_1,\qqq y_0\in \mH_{j+1}, \qqq y\in \mH_{j+2},
\\
-{1\/\r^2}(\r^2y_0')'+V y_0=0,\qqq
\\
-{1\/\r^2}(\r^2y_1')'+V y_1=1,\qqq y_1(0)=y_1(1)=0.
\end{aligned}
\]
Arguing as above, we have $y_0=0$. We consider $y_1$, which has the form
$$
y_1(x)=\int_0^1R(x,t)dt>0,
$$
$R(x,t)$ being the Green function for the last equation in \er{eB1}. Note that
 $R(x,t)\ge 0$ for any $x,t\in [0,1]\ts
[0,1]$, since the first eigenvalue of the Sturm-Liouville problem
equation
\[
\lb{SLeq}-{1\/\r^2}(\r^2\p')'+V\p=\l \p, \qqq \p(0)=\p(1)=0,
\]
is positive. We have $y_1'(0)\ge 0$ and  $y_1'(1)\le 0$, which
together with the equality
\[
\lb{y10} \int_0^1y_1''dx=y_1'(1)-y_1'(0)=0
\]
yields $y_1'(0)=y_1'(1)=0$. Let $\vp (x,t)$ be the solution of the
equation
$$
-\vp''-2\a q(x+t)\vp'+V(x+t)\vp =\l \vp ,\qqq \vp (0,t)=0, \qq \
\vp' (0,t)=1,\qq  t\in [0,1].
$$
Since $y_1(0)=y_1'(0)=0$, the solution of the last equation in
\er{eB1} has the form
$$
y(x)=-\int_0^x\vp (x-t,t)\cdot 1dt<0,
$$
which implies $y_1=0$, since $\vp (x,t)>0, x,t\in [0,1]$. $\BBox$

\subsection {Proof of Theorem \ref{T3} i)} In order to prove Theorem \ref{T3} i), we check all
conditions A)-D) in Theorem \ref{TA97} for the mapping $q\to P(q)$.

The statements A)  and B) have been proved in Lemma \ref{TAn}.

C) The two-sided  estimates \er{t3e} were proved in Lemma
\ref{Tpq1}.

 D) Let $q^\n\to q$ weakly in $\mH_j$ as $\n\to \iy$. Then we
deduce that $q^\n\to q$ strongly in $\mH_{j-1}$ as $\n\to \iy$ since
the imbedding mappings $\mH_j\to \mH_{j-1}$ are compact. Hence the
mapping $q\to P(q)-q'$ is compact.

Hence all conditions in Theorem \ref{TA97} hold true and the mapping
 $P:\mH_j \to \mH_{j-1}$ is a real analytic isomorphism  between
 the Hilbert spaces $\mH_j$ and $\mH_{j-1}$.
\BBox


\section {The mapping $P$ for singular case}
\setcounter{equation}{0}

 \subsection {Proof of Theorem \ref{T3} ii)}
 We consider the mapping $P$ given by
\[
\begin{aligned}
\lb{dPx}
& z'=P(q)=q'+q^2+h^2-c_0,\\
& h=h_0e^{-{2\/m}Q},\qq h_0\ge0, \qq Q(x)=\int_0^xq(t)dt,\qq
 c_0=\int_0^1(q^2+h^2)dx,
\end{aligned}
\]
for the singular case $q\in \mH_0$, where it is convenient to use
$z\in \mH_0$.

Consider the case $m\ge 1$. By Lemma \ref{TAn}, the mapping $P:
\mH_j\to \mH_{j-1}, j\ge 0$,  is a real analytic  and locally
invertible. We show that  $P$ is an injection.

Assume that $P$ is not injective. Then there exist $q, \wt q\in
\mH_0, q\neq \wt q$, such that $P(q)=P(\wt q)$. Then we have
$P(W)=\cB_\ve(z)=P(\wt W)$ where $\cB_\ve(z)=\{e\in \mH_0:
\|z-e\|<\ve \}$ for some domains $W, \wt W$ satisfying
\HOX{Corrected 2019.01.28} $q\in W\ss \mH_0$, $\wt q\in\wt
W\ss \mH_0$ and $W\cap \wt W=\es$, since $P$ is a real analytic  and
locally invertible. If we take any $q_1\in W\cap \mH_1$, then
$P(q_1)\in \cB_\ve(z)$ and there exists $\wt q_1\in \wt W\cap \mH_1$
such that $P(q_1)=P(\wt q_1), q_1\ne \wt q_1$. This yields
contradiction since $P: \mH_1\to \mH_0$ is a bijection.
\HOX{Corrected 2019.01.18} We have thus proven that  $P$ is
an injection.

We show the estimate \er{pq01}. We prove that of $z$ in terms of
$q$. Rewrite \er{dPx} in the form
\[
\lb{dPxy} z=q+\int_0^x\lt(f(t)-\int_0^1f(s)ds\rt)dt+
\int_0^1\lt(t-{1\/2}\rt)f(t)dt,
\]
where $f=q^2+u$. We further rewrite $z$ in the
form
\[
\lb{dPxya}  z=q+F-c_*,\qqq c_*=\int_0^1\lt(t-{1\/2}\rt)f(t)dt,
\]
 and
$$
\begin{aligned}
F(x)=\int_0^x(f(t)-\int_0^1f(s)ds)dt=\int_0^x(q^2(t)-\|q\|^2)dt+F_0(x),
\\
F_0(x)=\int_0^x(u(t)-\int_0^1u(s)ds)dt=A\int_0^1\int_0^x(e^{-\b
Q_t}-e^{-\b Q_s})dtds.
\end{aligned}
$$
Note here the following identity and estimate for some $\t\in [t,s]$:
$$
\begin{aligned}
e^{-\b Q_t}-e^{-\b Q_s}=\b (Q_s-Q_t)e^{-\b Q_\t},
\\
|e^{-\b Q_t}-e^{-\b Q_s}|\le \b \|q\|e^{\b \|q\|}.
\end{aligned}
$$
Then combining these estimates we obtain
\[
\lb{} |F(x)|\le 2\|q\|^2+2\|q\|+\b A\|q\|e^{\b \|q\|}.
\]
 Then we get
\[
\lb{pqz}
\begin{aligned}
\|z\|^2=(z,z)=(q+F-c_*,z)=(q+F,z)\le \|q\|\|z\|+(F,z)
\\
\le \|q\|\|z\|+ \|z\|\|q\|(2\|q\|+2+\b Ae^{\b \|q\|}),
\end{aligned}
\]
which yields \er{pq01}.

Let $m=1$. We prove that the mapping $P$ is a surjection.

We prove \er{t4e}. We show by an explicit construction that
for each $z\in \mH_0$, the equation \er{dPx}, i.e., $
z'=P(q)=q'+q^2+u(Q)-c_0(q)$, has a solution $q\in \mH_0$. For fixed
$z\in \mH_0$ we consider the auxiliary equation
\[
\lb{fp1} -\f''+z' \f=\l \f,\ \ \ x\in\R.
\]
We need the following results from \cite{K03}: there exists a unique
solution $\f(x)$ of \er{fp1} with unique $\l=\l_0(z)<\l_0^+(z)$
having the form
\[
\lb{fp2} \f(x)=e^{x}\f_1(x),\ \ \ \f_1(x+1)=\f_1(x),\ \ \
\phi_1(x)>0,\ \ \ \f_1'\in L^2(\T), \ \ \ \|\f_1\|=1.
\]
Here  $\l_0^+(z)$ is the lowest eigenvalue of the periodic problem
for \er{fp1}. Note that the function $\l_0(\cdot)$ is continuous on
$\mH_0$.

Let $z\in \mH_0$. Then we have $\f$. We show that there exists $q\in
\mH_0$ such that
\[
\lb{fp3} q+{h_0\/v}={\f'\/\f},\ \ \ {\rm where}\ \ \
v(x)=e^{2\int_0^xq(s)ds}.
\]
The equation \er{fp3} is  linear   for  $v$, since
$v'=2qv$ and
\[
\lb{fp4} {v'\/2v}+{h_0\/v}={\f'\/\f} \ \ \ \Leftrightarrow\ \ \
v'+2h_0=2{\f'\/\f}v.
\]
This equation has a periodic solution of the form
\[
\lb{fp5} v(x)=C_1\int_0^x\lt({\f(x)\/\f(t)}\rt)^2dt
+C_2\int_x^1\lt({\f(x)\/\f(t)}\rt)^{2}dt  ,\ \ \ x\in\R,
\]
where $$
C_1={2h_0\/e^{2}-1},\qqq C_2=C_1e^2.
$$
Thus we have $v'\in L^2(0,1)$ and $ v>0$. Then we obtain
$2q={v'\/v}\in \mH_0$ and $h=h_0e^{-2\int_0^xq(s)ds}$.

We show that $z,q$ satisfy the equation \er{dPx}. Differentiating
\er{fp3}  and using $h=h_0/v$ we obtain
\[
\lb{fp6}
\begin{aligned}
& (q+h)'={\f''\/\f}-(q+h)^2=(p'-\l_0)-(q+h)^2,\\
& q'-2qh=(z'-\l_0)-(q^2+2qh+h^2),\\
& z'=q'+q^2+h^2+\l_0,
\end{aligned}
\]
where
\[
\lb{fp7} -\l_0(z)=\int_0^1(q^2+h^2)dx=\|q\|^2+\|h\|^2.
\]
Thus
\HOX{Corrected 2019.01.28} $q$ satisfies equation \er{dPx} and  $P$ is a
surjection. Hence $P:\mH_0\to \mH_{-1}$ is a real analytic
isomorphism between $\mH_0$ and $\mH_{-1}$.

Now we prove the estimate (\ref{pq01}). We show that it is
sufficient to consider $z\in \mH_1$. Substituting $\f=e^{x}\f_1$
into \er{fp1}, we obtain the equation for $\f_1$:
\[
\label{f1} -\phi_1''-2\phi_1'+z'\phi_1=(\l_0(z)+1)\phi_1,\ \ \
(x,z)\in\R\ts \mH_1.
\]
Recall $\|\phi_1\|=1$.  Multiplying \er{f1} in $L^2(0,1)$ by
$\phi_1$ we obtain the identity
$$
\|\phi_1'\|^2+(z'\phi_1,\phi_1)=\l_0(z)+1.
$$
Repeating the arguments from \cite{K03} we obtain the estimate
\[
\label{ff1}
 -\l_0(z)-1\le 2\|z\|^2(1+2\|z\|^2),\ \ \ z\in\mH_1.
\]
By the continuity of $\l_0(\cdot)$ on $\mH_0$  the estimate \er{ff1}
holds true for $z\in \mH_0$.  Substituting \er{fp7} into \er{ff1} we
get
\[
\label{2.15.13} \|q\|^2+\|h\|^2 -1  \le 2\|z\|^2(1+2\|z\|^2) .
\]
Equation (\ref{fp4}) and $\phi=e^{x}\phi_1$ give
$$
 \int_0^1h(x)dx=\int_0^1{\phi'(x)\/\phi(x)}dx={1},
$$
which implies $\|h\|^2\ge 1$. Then (\ref{2.15.13}) gives $\|q\|^2\le
2\|z\|^2(1+2\|z\|^2),$ which yields  \er{pq01}. \BBox



{\bf Proof of Theorem \ref{T2}.} We recall that due to  Theorem
\ref{Tgap}  in \cite{K03}  the mapping
\[
\lb{FKCx} \P: p\mapsto \P(p)= (\P_n(p))_{n=1}^{\iy}
\]
is a real-analytic isomorphism between $\mH_{-1}$ and
$\ell_{-1}^2\os \ell_{-1}^2$.

i) Fix $m\ge 1$. Then due to Theorem \ref{T3}, the mapping
$P:\mH_0\to \mH_{-1}$, given by \er{dP1}, is real analytic and
injective. Thus due to \er{gm}--\er{f4} and Theorem \ref{Tgap} and
the relation between $q$ and $P$ given by \er{a7x}, \er{6}, we
obtain the identity
\[
\lb{aa1} \p(q)=\P(P(q)),\qq \forall \ q\in \mH_0.
\]
The mapping $\p(\cdot)$ is the composition of two mappings $\P$ and
$P$. The properties described above
\HOX{Corrected 2019.01.28} show that the mapping
$\p:\mH_{0}\to \ell_{-1}^2\os \ell_{-1}^2$ is real analytic and
injective. Moreover, combining \er{pq01} and \er{1.7}, we  obtain
\er{13x}.

ii) Let $m=1$. Recall that  Theorem \ref{T3} (in our case $j=0$)
\HOX{Corrected 2019.01.18} proves that the mapping $q\to P=P_\n(q)$ defined by
\er{6} is a real-analytic isomorphism between $\mH_0$ and
$\mH_{-1}$. Thus by virtue of \er{gm}--\er{f4}, Theorem \ref{Tgap}
and the relation between $q$ and $P$ given by \er{a7x}, \er{6}, we
obtain the following equality
$$
\p(q)=\P(P(q)),\qq \forall \ q\in \mH_0.
$$
The mapping $\p(\cdot)$ is the composition of two mappings $\P$ and
$P$,
\HOX{Corrected 2019.01.28} each of which  is a corresponding analytic isomorphism
(see Theorem \ref{Tgap} and Theorem \ref{T3}). Then the mapping $
\p:q\mapsto (\P_n(P(q))_{n=1}^{\iy}, $ is a real-analytic
isomorphism between $\mH_0$ and $\ell_{-1}^2\os \ell_{-1}^2$.
Finally, combining \er{t4e} and \er{1.7}, we  obtain \er{14x}.
 \BBox

{\bf Proof of Corollary \ref{Todd}.} i) It is clear that the mapping
$P$ given by  \er{dP1} satisfies $P:\mH_j^{odd}\to
\mH_{j-1}^{even}, j\ge1$. Arguments in the proof of Theorem
\ref{T3} prove that this  mapping $P$
 is a real analytic
isomorphism between the Hilbert spaces $\mH_j^{odd}$ and
$\mH_{j-1}^{even}$.

ii) In order to prove ii) we need the Garnett-Trubowitz' results
\cite{GT84} or \cite{GT87} on the gap-length mapping for
Schr\"odinger operator $Sy=-y''+py$ on the torus $\T$ with the
potential $p\in \mH_0^{even}$. Define the mapping $\P^e:\mH_0^{even}
\to \ell^2$ by
\[
\lb{GTx} p\to \P^e(p)=(\P_n^e(p))_1^\iy,\qqq
\P_n^e(p)={\l_n^+(p)+\l_n^-(p)\/2}-\m_n(p),
\]
for all $n\ge 1$. Here $\l_n^\pm(p)$ is the periodic and
antiperiodic eigenvalues for the equation $y''+py=\l y$ and
$\m_n(p)\in [\l_n^+(p),\l_n^-(p)]$ are the corresponding Dirichlet
eigenvalues for the problem $y''+py=\l y, y(0)=y(1)=0$. By
Garnett-Trubowitz \cite{GT84}, we then see that the mapping
$\P^e:\mH_0^{even}\to\ell^2$ given by \er{GTx}
 is a
\HOX{Corrected 2019.01.28} real analytic isomorphism between the Hilbert spaces
$\mH_1^{odd}$ and $\ell^2$. Repeating the argument in the proof of
Theorem \ref{T1}, we prove Corollary \ref{Todd} ii) for the case
$\mH_1^{odd}$.

\HOX{Corrected 2019.01.28}
We consider the singular case $\mH_0^{odd}$ and deal with only the case $m=1$. It is
clear that the mapping $P$ given by \er{dP1} satisfies
$P:\mH_0^{odd}\to \mH_{-1}^{even}$. The results from Theorem
\ref{T3} and arguments in the proof of Theorem \ref{T3} give that
this mapping $P$
 is a real analytic
isomorphism between the Hilbert spaces $\mH_0^{odd}$ and
$\mH_{-1}^{even}$.

From  Theorem \ref{Tgap} we deduce that the mapping
$\P^e:\mH_{-1}^{even}\to\ell_{-1}^2$ given by \er{GTx}, is a
real analytic isomorphism between the Hilbert spaces
$\mH_{-1}^{odd}$ and $\ell_{-1}^2$. Thus
\HOX{Corrected 2019.01.28}
by the composition of the
mapping $P:\mH_0^{odd}\to \mH_{-1}^{even}$ and
$\P^e:\mH_{-1}^{even}\to\ell_{-1}^2$, the mapping
$\p^e:\mH_0^{odd}\to\ell_{-1}^2$, given by \er{epj}
 is a  real analytic isomorphism between the Hilbert spaces
$\mH_0^{odd}$ and $\ell_{-1}^2$.  \BBox


\subsection {Proof of Corollary  \ref{T4}}
  Due to Theorem \ref{T3} the mapping $v\to \cE_1$  is a real
analytic isomorphism between the spaces $\mH_j$ and $\mH_{j-1}$.
Moreover, for each $\cE_1\in \mH_0$  there exists a unique  $v\in
\mH_1$ and the constant $\cE_{0}$ in \er{Rcr1}, which   is uniquely
defined by $\cE_1$.  We show estimates. Let $\cE_1=-v'-g$, where
$g={v^2\/m}-\cE_0$ and $\cE_0=-{\|v\|^2\/m}$. We have
\[
\begin{aligned}
\lb{es1}
&\|\cE_1\|^2=\|v'+g\|^2=\|v'\|^2+\|g\|^2+2(v',g)\\
& \ \ \ =\|v'\|^2+\|g\|^2=\|v'\|^2+{1\/m^2}\|v^2\|^2-\cE_0^2,
\end{aligned}
\]
since
$$
\begin{aligned}
 &(v',g)={1\/m}(v',v^2)={1\/m}\int_0^1v'v^2dx=0,
\\
&\|g\|^2={1\/m^2}\|v^2\|^2+\cE_0^2-{2\/m}\|v\|^2\cE_0={1\/m^2}\|v^2\|^2-\cE_0^2.
\end{aligned}
$$
This identity \er{es1} and the following estimate
$$
\|v^2\|^2=\int_0^1v^4dx\le \|v\|^2\sup_{x\in \T}|v(x)|^2\le
\|v\|^2\|v'\|^2,
$$
 yield
\[
\lb{es2} \|v'\|^2\le \|\cE_1\|^2\le
\|v'\|^2+{1\/m^2}\|v\|^2\|v'\|^2-\cE_0^2.
\]
\BBox


\bigskip

\setlength{\itemsep}{-\parskip} \footnotesize

\no  {\bf Acknowledgments.} We thank Andrei Badanin for the Fig.
\er{Fig1}. Various parts of this paper were written during Evgeny
Korotyaev's stay  in the Mathematical Institute of University of
Tsukuba, Japan. He is grateful to the institute for the hospitality.
  H. Isozaki is supported by Grants-in-Aid for
Scientific Research (S) 15H05740, and (B) 16H03944, Japan Society
for the Promotion of Science. E. Korotyaev  is supported  by the RSF
grant No. 18-11-00032.


\begin{thebibliography}
{99999}\setlength{\itemsep}{-\parskip} \footnotesize


\bibitem{AA07} Aberra, D.; Agrawal, K. Surfaces of revolution in n
dimensions. Internat. J. Math. Ed. Sci. Tech.  38  (2007),  no. 6,
843-851.



\bibitem {BKK03}  Badanin, A.; Klein, M.; Korotyaev, E.
The Marchenko-Ostrovski mapping and the trace formula for the
Camassa-Holm equation. J. Funct. Anal. 203 (2003), no. 2, 494--518.

\bibitem{B94}
Bleher, P. M. Distribution of energy levels of a quantum
       free particle on a surface of revolution. Duke Math. J.  74  (1994),
       no. 1, 45�-93.


\bibitem {BH84} Br\"uning, J.; Heintze, E. Spektrale Starrheit
  gewisser Drehfl\"achen, Math. Ann. 269 (1984) 95--101.







\bibitem{CY76}  Cheng, S.Y.; S.T. Yau,  On the regularity
of the solution of the n-dimensional Minkowski problem, Comm. Pure
Appl. Math. 29(1976), 495--516.


\bibitem{C06}  Cheshkova, M. A. Surfaces of revolution of
 constant Gaussian curvature.   (Russian) Differentsial'naya
 Geom. Mnogoobraz. Figur  No. 37  (2006), 176�179.

\bibitem{CM93}  Coleman C., McLaughlin J.: Solution of the inverse problem
for an impedance with integrable derivative, Commun. Pure Appl.
Math.  46(1993), 145--184; 185--212.



\bibitem{DK09}   Dorfmeister, J.; Kenmotsu, K. Rotational
hypersurfaces of periodic mean curvature. Differential Geom. Appl.
27  (2009),  no. 6, 702�-712.

\bibitem{DK08}    Dorfmeister, J.; Kenmotsu, K. On a theorem
 by Hsiang and Yu. Ann. Global Anal. Geom.  33  (2008),  no. 3, 245�-252.

\bibitem {D75} Dubrovin, B. Periodic problems for the Korteveg de Vries
 equation in the class of  finite-gap potentials.
Funct. Anal. Appl. 9(1975), 215--223.


\bibitem{E98}
Engman, M. Sharp bounds for eigenvalues
  and multiplicities on surfaces of revolution. Pacific J. Math.  186(1998),
   no. 1, 29�-37.



\bibitem{GT84} Garnett J.; Trubowitz E. Gaps and bands of one dimensional
 periodic Schr\"odinger operators. Comment. Math. Helv. 59(1984), 258--312.

\bibitem{GT87} Garnett J.; Trubowitz E.  Gaps and bands of one dimensional
 periodic Schr\"odinger operators II. Comment. Math. Helv. 62 (1987), 18--37.












\bibitem{IM75} A. Its, V. Matveev: Schr\"odinger operator with the finite-gap
 spectrum and the N-soliton solutions of the Korteveg de Fries
 equation. Theoret. Math. Phys. 23 (1975), 343--355.

\bibitem{IK17} Isozaki, H. ;  Korotyaev E. Inverse spectral theory and
the Minkowski problem for the surface of revolutions, Dynamics of
PDE, 14 (2017),  No. 4, 321--341.

\bibitem{IK17x} Isozaki, H. ;  Korotyaev E. Global transformations preserving
Sturm-Liouville spectral data, Rus. J. Math. Phys.
 24 (2017), no. 1, 51-68.

\bibitem{KK97} Kargaev, P.;  Korotyaev, E. Inverse problem for the Hill
operator, the direct approach.  Invent. Math., 129 (1997), no. 3,
567--593.

\bibitem{KK99}
Kargaev, P.;  Korotyaev, E. Erratum The inverse problem for the Hill operator,
 a direct approach, Invent. math. 138(1999), 227.

\bibitem{Ke03}   Kenmotsu, Katsuei Surfaces of revolution with
periodic mean curvature. Osaka J. Math.  40  (2003),  no. 3,
687-696.

\bibitem{KK00}
Klein M., Korotyaev E.: Parametrization of periodic weighted
operators in terms of gaps lengths. Inverse Problems, 16(2000),
no.6, 1839-1860.

\bibitem{KT96}
  Konopelchenko, B. G.; Taimanov, I. A. Constant mean curvature surfaces
  via an integrable dynamical system. J. Phys. A  29  (1996),  no. 6, 1261�-1265.

\bibitem{K97}   E. Korotyaev:
The inverse problem for the Hill operator I.
  Inter. Math. Research. Notes. 3 (1997), 113-125.

\bibitem{K98}
Korotyaev, E. Estimates of periodic potentials in terms of gap
lengths. Comm. Math. Phys. 197 (1998), no. 3, 521--526.

\bibitem{K99} Korotyaev, E. Inverse problem and the trace formula for
the Hill operator. II. Math. Z. 231 (1999), no. 2, 345--368.

\bibitem{K00e} Korotyaev, E. Estimates for the Hill operator. I. J.
Differential Equations 162 (2000), no. 1, 1--26.

\bibitem{K00}  Korotyaev, E. Inverse Problem for Periodic "Weighted" Operators.
J. Funct. Anal. 170(2000), no. 1, 188--218.

\bibitem{K02}   Korotyaev, E. Invariance principle for inverse problems.
Int. Math. Res. Not. 2002, no. 38, 2007--2020.

\bibitem{K03}  Korotyaev, E.
Characterization of the spectrum of   Schr\"odinger operators with
periodic distributions.  Int. Math. Res. Not. 2003, no. 37,
2019--2031.


\bibitem{K03p} Korotyaev, E. Periodic "weighted" operators. J.
Differential Equations 189 (2003), no. 2, 461--486.

\bibitem{Kr57} Krein M. On the characteristic function
$A(\l)$ of a linear canonical system differential equation of the
second order with periodic coefficients (Russian),  Prikl. Mat. Meh.
21(1957), 320--329.


\bibitem{L87} Levitan, B. Inverse Sturm-Liouville problems.
Utrecht: VNU Science Press, 1987.

\bibitem{M86} Marchenko, V. Sturm-Liouville operator and
applications. Basel: Birkh\"auser, 1986.


\bibitem{MO75}
Marchenko, V. A.; Ostrovski, I. V.
A characterization of the spectrum of the Hill operator. (Russian)
Mat. Sb. (N.S.) 97(139) (1975), no. 4(8), 540--606, 633--634.

\bibitem{M81}  J. Moser,  Integrable Hamiltonian system and spectral theory.
Academia Nationale dei Liecei, Scuola Superiore, Pisa 1981.

\bibitem{N74}
Novikov, S. The periodic problem for the Korteveg de Fries
equation. Funct. Anal. Appl. 8(1974), 236-246.


\bibitem{P74}
A. V. Pogorelov,
The Minkowski multidimensional problem. (Russian) [Hilbert's fourth
problem] Izdat. ``Nauka'', Moscow, 1974. 79 pp.
   English translation by V. Oliker. Introduction by Louis Nirenberg.
   Scripta Series in Mathematics. V. H. Winston $\&$
    Sons, Washington, D.C.; Halsted Press [John Wiley $\&$ Sons],
    New York-Toronto-London, 1978.


\bibitem{PT87}
P\"oschel, J.; Trubowitz, E. Inverse spectral theory.
Pure and Applied Mathematics, 130. Academic Press, Inc., Boston, MA,
1987.

\bibitem{T77}
Trubowitz, E. The inverse problem for  periodic potentials. Commun.
Pure Appl. Math. 30(1977), 321-337.


\bibitem{SW03}
Sanders, J.; Wang, J. Integrable systems in n-dimensional
Riemannian geometry. Mosc. Math. J.  3  (2003),  no. 4, 1369�-1393.






\bibitem{Ta97}
Taimanov, I. Surfaces of revolution in terms
of solitons. Ann. Global Anal. Geom.  15  (1997),  no. 5, 419-�435.



\bibitem{Z98}
Zelditch, S. The inverse spectral problem for surfaces
of revolution. J. Differential Geom.  49  (1998),  no. 2, 207�-264.




\end{thebibliography}
\end{document}